\documentclass[11pt,reqno]{amsproc}

\title[]{Nernst-Planck-Navier-Stokes Systems Far From  Equilibrium}

\author{Peter Constantin}
\address{Department of Mathematics, Princeton University, Princeton, NJ 08544}
\email{const@math.princeton.edu}

\author{Mihaela Ignatova}
\address{Department of Mathematics, Temple University, Philadelphia, PA 19122}
\email{ignatova@temple.edu}

\author{Fizay-Noah Lee}
\address{Program in Applied and Computational Mathematics, Princeton University, Princeton, NJ 08544}
\email{fizaynoah@princeton.edu}

\usepackage[margin=1in]{geometry}
\usepackage{amsmath, amsthm, amssymb}
\usepackage{times}
\usepackage{color}
\usepackage{hyperref}

\newcommand{\pa}{\partial}
\newcommand{\la}{\label}
\newcommand{\fr}{\frac}
\newcommand{\na}{\nabla}
\newcommand{\be}{\begin{equation}}
\newcommand{\ee}{\end{equation}}
\newcommand{\ba}{\begin{array}{l}}
\newcommand{\ea}{\end{array}}
\newcommand{\Rr}{{\mathbb R}}

\newcommand{\beg}{\begin}

\renewcommand{\div}{{\mbox{div}\,}}
\newcommand{\D}{\Delta}
\newcommand{\C}{C_{\Gamma}}
\date{today}
\begin{document}
\begin{abstract}
We consider ionic electrodiffusion in fluids, described by the Nernst-Planck-Navier-Stokes system. We prove that the system has global smooth solutions for arbitrary smooth data: aribitrary positive Dirichlet boundary conditions for the ionic concentrations,  arbitrary Dirichlet boundary  conditions for the potential, arbitrary positive initial concentrations, and arbitrary regular divergence-free initial velocities. The result holds for any positive diffusivities of ions, in bounded domains with smooth boundary in three space dimensions, in the case of two ionic species, coupled to Stokes equations for the fluid. The result also holds in the case of Navier-Stokes coupling, if the velocity is regular.  The global smoothness of solutions is also true for arbitrarily many ionic species, if all their diffusivities are the same. 
\end{abstract}

\noindent

\maketitle

\section{Introduction}
The Nernst-Planck-Navier-Stokes system describes the evolution of ions in a Newtonian fluid \cite{rubibook}. Several species of ions, with different valences $z_i$ diffuse with diffusivities $D_i>0$, and are carried by an incompressible fluid with constant density and with velocity $u$, and by an electrical field generated by the local charge $\rho$ and by voltage applied at the boundaries. The system of equations is
\be
\pa_t c_i + u\cdot\na c_i = D_i \div (\na c_i + z_ic_i\na \Phi),
\la{cieq}
\ee
$i=1,2, \dots m$, coupled to the Poisson equation
\be
-\epsilon\D\Phi = \sum_{i=1}^m z_i c_i = \rho,
\la{poi}
\ee
and to  the Navier-Stokes equations
\be
\pa_t u + u\cdot\na u -\nu\D u + \na p = -K \rho\na\Phi, \quad \na\cdot u = 0,
\la{nse}
\ee
or to the Stokes equation,
\be
\pa_t u -\nu\D u + \na p = -K \rho\na\Phi, \quad \na\cdot u = 0.
\la{se}
\ee
The function $c_i(x,t)$ represents the local concentration of the $i$-th species, and $\Phi$ is an electrical potential created by the charge density $\rho$. The positive constant $\epsilon$ is proportional to the square of the Debye length.  The kinematic viscosity of the fluid is $\nu>0$ and
$K>0$ is a coupling constant with units of energy per unit mass. This constant is proportional to the product of the Boltzmann's constant $K_B$ and the absolute temperature $T_K$. The potential $\Phi$ has been nondimensionalized so that $\fr{k_BT_K}{e}\Phi$ is the physical electrical potential, where $e$ is elementary charge. The charge density $\rho$ has been nondimensionalized so that $e\rho$ is the physical electrical charge density.

The boundary conditions for $c_i$ are inhomogeneous Dirichlet,
\be
{c_i(x,t)}_{\left | \right. \pa\Omega} = \gamma_i(x)
\la{gammai}
\ee
The boundary conditions for $\Phi$ are inhomogeneous Dirichlet
\be
\Phi(x,t)_{\left | \right. \pa\Omega} = W(x)
\la{phibc}
\ee
and the boundary conditions for the Navier-Stokes or Stokes equations are homogeneous Dirichlet,
\be
u_{\left | \right. \pa\Omega} = 0.
\la{ubc}
\ee
The bounded connected domain $\Omega$ need not be simply connected. The functions $\gamma_i$ defined on the boundary of the domain are given, smooth positive time independent functions. We denote by $\Gamma_i$ positive, smooth, time independent extensions of these functions in the interior 
\be
\gamma_i = {\Gamma_i}_{\left | \right. \pa\Omega} 
\la{Gammai}
\ee
The function $W$ is also a given and time independent smooth function. 

The NPNS system is a well posed semilinear parabolic system. The question we are discussing is whether it has global smooth solutions, that is, whether given
arbitrary smooth initial data and arbitrary smooth boundary conditions, do smooth solutions exist for all time, or do some solutions blow up.

As it is well known, semilinear parabolic scalar equations can blow up in finite time. The simplest such example is a semilinear heat equation (\cite{gigak}). Semilinear systems in which there is a single concentration carried by the gradient of a potential can also blow up. Well-known examples are chemotaxis equations, such as the Keller-Segal equation (\cite{perthame}).  The system we are discussing involves the Navier-Stokes equation,
where the question we discuss here is a major open problem, but even in the absence of fluid or in the case of Stokes flow coupled to Nernst-Planck equations, the problem of global existence of smooth solutions discussed here remains open in its full generality.

Boundary conditions play an essential role in the behavior the solutions of 
the NPNS system. No flux (blocking) boundary conditions for the concentrations model situations in which the boundaries are impermeable to the ions.  Dirichlet (selective) boundary conditions for the concentrations discussed in this paper model situations in which boundaries maintain a certain concentration of ions.   Global existence and stability of solutions of the Nernst-Planck equations, uncoupled to fluids has been obtained in \cite{biler}, \cite{choi}, \cite{gajewski} for blocking boundary conditions in two dimensions, or in three dimensions for small data, or in a weak sense. 
The system coupled to fluid equations was studied in  {\cite{schmuck}} where  global existence of weak solutions is shown in two and three dimensions for homogeneous Neumann boundary conditions on the potential, a situation without boundary current. In {\cite{ryham}}, homogeneous Dirichlet boundary conditions on the potential are considered, and global existence of weak solutions is shown in two dimensions for large initial data and in three dimensions for small initial data (small perturbations and small initial charge). In \cite {bothe} the problem of global regularity in two dimensions was considered for Robin boundary conditions for the potential. In \cite{ci} global existence of smooth solutions for blocking boundary conditions and for uniform selective (special, stable Dirichlet) boundary conditions were obtained in two space dimensions. In \cite{np3d} blocking and uniform selective boundary conditions were used to prove nonlinear stability of Boltzmann states in three space dimensions.

While blocking and uniformly selective boundary conditions lead to stable configurations, instabilities may occur for general selective boundary conditions. These instabilities have been studied in simplified models mathematically and numerically (\cite{rubizaltz}, \cite{zaltzrubi}) and observed in physical experiments \cite{rubinstein}.  In this paper we consider the system with large data, with general selective boundary conditions, in situations in which instabilities may occur. We prove global regularity of solutions for two cases: if there are only two species (cations and anions, $m=2$) or if there are many species, but they all have the same diffusivities ($D_1= \dots = D_m$). The difficulty in three dimensions, even when there is no fluid, is in bounding the nonlinear growth of the concentrations.

This paper is organized as follows. In Section \ref{prel}  we prove a necessary and sufficient condition for global regularity of solutions. In the case of the Nernst-Planck system coupled to Stokes equations this condition (Theorem \ref{bkm}) states that, if (and only if) 
\be
\sup_{0\le \tau <T}\int_0^\tau\left(\int_{\Omega}|\rho(x,t)|^2dx\right)^2dt  = B(T) <\infty
\la{regcond}
\ee
is finite, then the solution is smooth on $[0,T]$.  
In the case of coupling to the Navier Stokes equation, the condition (\ref{regcond}) is supplemented by a well-known condition for regularity for the Navier-Stokes equations. 

In Section \ref{en} we introduce functionals of $(c_i, \Phi)$ which are used to cancel the contribution of electrical forces in the Navier-Stokes or Stokes energy balance, at the price of certain quadratic error terms. The first functional is a sum of relative entropies and a potential norm. The second one, which is just the potential part of the first one, is weaker and has a weaker dissipation, but it introduces better error terms.  

Section \ref{glob} is devoted to quadratic bounds, which imply global regularity by the criterion established in Theorem \ref{bkm}. It is only here that the restriction to $m=2$ (Theorem \ref{twospecies}) or to $D_1 = \dots = D_m$ (Theorem \ref{equaldiffs}) is used. In these two special circumstances we  show that there is a cubic dissipation term proportional to $\|\rho\|_{L^3}^3$ in the evolution of the quadratic norms that we are employing to control the concentrations. This cubic term is essential, because although it has a small prefactor, it can be used to absorb all the quadratic errors ocurring in the second energy.

\section{Preliminaries}\la{prel}
We denote by $C$ absolute constants. We denote $C_{\Gamma}$ any constant that depends only on  the parameters and boundary data of the problem, i.e. on $\nu, K, \epsilon, z_i, D_i,$, on the domain $\Omega$ itself, on norms of $W$ and on norms of $\Gamma_i$. These constants may change form line to line, and they are explicitly computable. They do not depend on solutions, or initial data. We do not keep track of them to ease notation and focus on the ideas of the proofs. 

We consider a bounded domain $\Omega\subset\Rr^3$ with smooth boundary. We denote space $L^p(\Omega) = L^p$  and norms simply by $\|\cdot\|_{L^p}$.  We denote by $H = L^2(\Omega)^3\cap \{u\left |\right. \; \div u = 0\}$ the space of square integrable, divergence-free velocities with norm $\|\cdot\|_H$ and by $V = H_0^1(\Omega)^3\cap \{ u\left | \right.\; \div u =0\}$ the space of divergence free vectors fields with components in $H_0^1(\Omega)$, with norm $\|\cdot\|_V$. We denote by $\mathbb P$ the Leray projector $\mathbb P : L^2(\Omega)^3 \to H$, and by $A$ the Stokes operator
\be
A = -\mathbb P\D, \quad \quad A: \mathcal{D} (A)\to H
\la{Aop}
\ee
where 
\be
\mathcal D(A) = H^2(\Omega)^3\cap V.
\la{mathcalda}
\ee

\beg{defi} \la{strongsol}We say that $(c_i, \Phi, u)$ is a strong solution of the system
(\ref{cieq}), (\ref{poi}), (\ref{nse}) or (\ref{cieq}), (\ref{poi}), (\ref{se}) with boundary conditions (\ref{gammai}), (\ref{phibc}), (\ref{ubc}) on the time interval $[0,T]$ if 
$u\in L^{\infty}(0,T; V)\cap L^2(0,T; \mathcal D(A))$ and $c_i\in L^{\infty}(0,T; H^1(\Omega))\cap L^2(0,T; H^2(\Omega))$ solve the equations in distribution sense and the boundary conditions in trace sense.
\end{defi} 

It is well-known that strong solutions of Navier-Stokes equations are as smooth as the data permit (\cite{cf}). The same is true for the Nernst-Planck-Navier-Stokes equations.

\beg{thm}\la{locex} Let $c_i(0)-\Gamma_i\in H_0^1(\Omega)$, and $u(0)\in V$. There exists $T_0$ depending on $\|u_0\|_V$ and $\|c_i(0)\|_{H^1(\Omega)}$, the boundary conditions $\gamma_i, W$, and the parameters of the problem ($\nu, K, D_i, \epsilon, z_i)$, so that the  system (\ref{cieq}), (\ref{poi}), (\ref{nse}) (or the system (\ref{cieq}), (\ref{poi}), (\ref{se})) with boundary conditions (\ref{gammai}), (\ref{phibc}), (\ref{ubc})  has a unique strong solution $(c_i, \Phi, u)$ on the interval $[0,T_0]$.
\end{thm}
\noindent\beg{proof}
We write
\be
c_i = q_i + \Gamma_i.
\la{qi}
\ee
and note that the equations  (\ref{cieq}) can be written as
\be
\pa_t q_i + u\cdot\na q_i = D_i \div (\na q_i + z_iq_i\na \Phi) + F_i
\la{Fiqieq}
\ee
where
\be
F_i = -u\cdot\na \Gamma_i + D_i\div(\na \Gamma_i + z_i \Gamma_i\na\Phi).
\la{Fi}
\ee
The boundary conditions for $q_i$ are homogeneous Dirichlet,
\be
{q_i}_{\left | \right. \pa\Omega} = 0.
\la{qibc}
\ee
We sketch only the apriori bounds for the proof. The actual construction of solutions can be done via Galerkin approximations.
Taking the scalar product of (\ref{Fiqieq}) with $-\D q_i$, we estimate the terms
\be
\ba
\left|z_iD_i\int_{\Omega}(\na q_i\cdot\na\Phi + q_i\D\Phi)\D q_i dx\right|
\le \C\left (\|\na\Phi\|_{L^6}\|\na q_i\|_{L^3} + \|q_i\|_{L^4}\|\rho\|_{L^4}\right)\|\D q_i\|_{L^2} \\
\le \C\left((\|\rho\|_{L^2} +1)\|\na q_i\|_{L^2}^{\fr{1}{2}}\|\D q_i\|_{L^2}^{\fr{3}{2}}  + \|q_i\|_{L^2}^{\fr{1}{4}} \|\rho\|_{L^2}^{\fr{1}{4}}\|\rho\|_{L^6}^{\fr{3}{4}}\|\na q_i\|_{L^2}^{\fr{3}{4}}\|\D q_i\|_{L^2}\right)
 \ea
\la{nonlqh1}
\ee
where we used
\be
\|\na\Phi\|_{L^6}\le \C(\|\rho\|_{L^2}+1)
\la{naphi6}
\ee
(the inhomogeneous boundary conditions are accounted for in the added term $\C$), embedding $H^1\subset L^6$ and interpolation. The advective term is estimated
\be
\left| \int_{\Omega} (u\cdot\na q_i)\D q_idx\right|\le C_\Gamma \|u\|_V\|\na q_i\|_{L^2}^\fr{1}{2}\|\D q_i\|_{L^2}^\fr{3}{2}
\la{conv}
\ee
The forcing term is estimated
\be
\left| \int_{\Omega}F_i \D q_idx\right |\le \C(\|u\|_H + \|\rho\|_{L^2} +1)\|\D q_i\|_{L^2}
\la{forcedel}
\ee
where we used
\be
\|\na\Phi\|_{L^p}\le \C(\|\rho\|_{L^2}+1)
\la{naphirholp}
\ee
valid for $p\le 6$.
Now we have
\be
\|\rho\|_{L^p}\le \C\sum_{i=1}^m(\|q_i\|_{L^p} + 1)
\la{normsrhoq}
\ee
and therefore, using also the Poincar\'{e} inequality we obtain
\be
\fr{d}{dt}\sum_{i=1}^m\|\na q_i\|^2_{L^2}  + \sum_{i=1}^m \|\D q_i\|_{L^2}^2 \le \C \left[ \left(\sum_{i=1}^m\|\na q_i\|^2_{L^2} + 1 \right)^3  + \|u\|_V^6\right].
\la{normqs}
\ee
We take the scalar
product of (\ref{nse}) with $Au$ we obtain, using well known estimates for the NSE \cite{cf},
\be
\fr{d}{dt} \|u\|_{V}^2 + \nu\|Au\|_{H}^2 \le \C(\|u\|_V^6 + \|\rho \na\Phi\|_{L^2}^2).
\la{nsevineq}
\ee
Adding to (\ref{normqs}) we obtain short time control of the norms 
required by the definition of strong solutions.
\end{proof}

\beg{lemma} \la{sufcondB}Let $(c_i, \Phi, u)$ be a strong solution of the system (\ref{cieq}), (\ref{poi}), (\ref{nse}) (or (\ref{cieq}), (\ref{poi}), (\ref{se})) with boundary conditions (\ref{gammai}), (\ref{phibc}), (\ref{ubc})  on the interval $[0,T]$. Let $2\le p <\infty$ be an even integer, take  $c_i(0)\in L^p$ and consider the quantity 
\be
\int_0^T \|\rho(t)\|_{L^2}^4 dt = B(T)<\infty,
\la{condB}
\ee
which is finite because $(c_i, \Phi, u)$ is a strong solution.
Then $c_i\in L^{\infty}(0,T; L^p)$ and 
\be
\sum_{i=1}^m \|c_i(t)\|_{L^p}\le \C\left(\sum_{i=1}^m \|c_i(0)\|_{L^p}+\int_0^T\|u\|_{L^p}dt +1\right)e^{\C\left(T +B(T)\right)}.
\la{lftylpb}
\ee
holds.
\end{lemma}
\beg{proof} We multiply the equation (\ref{Fiqieq}) by $q_i^{p-1}$ and integrate. We estimate the term
\be
\ba
\left|z_iD_i\int_{\Omega}\div(q_i\na\Phi)q_i^{p-1}dx\right| 
\le \C\|\na\Phi\|_{L^6}\|q_i^{\fr{p}{2}}\|_{L^3}\|\na q_i^{\fr{p}{2}}\|_{L^2}\\
\le \C(\|\rho\|_{L^2} +1)\|q_i^{\fr{p}{2}}\|_{L^2}^{\fr{1}{2}}\|\na q_i^{\fr{p}{2}}\|_{L^2}^{\fr{3}{2}}
\ea
\la{nonlqp}
\ee
where we used one integration by parts, allowed by the vanishing of $q_i$ at the boundary and interpolation. We estimate the forcing term
\be
\left| \int_{\Omega}F_iq_i^{p-1}dx\right | \le \C( \|u\|_{L^p} + \|\na\Phi\|_{L^p} + \|\rho\|_{L^p})\|q_i\|_{L^p}^{p-1}.
\la{fiqilp}
\ee
Using the dissipative term 
\be
D_i\int_{\Omega}\D q_i q_i^{p-1} dx = -D_i\fr{4(p-1)}{p^2}\int_{\Omega} |\na q_i^{\fr{p}{2}}|^2 dx,
\la{dispqlp}
\ee
and then discarding it, we obtain 
\be
\fr{d}{dt}\|q_i\|_{L^p}\le \C(\|\rho\|_{L^2}^4 +1)\|q_i\|_{L^p} + \C( \|u\|_{L^p} + \|\na\Phi\|_{L^p} + \|\rho\|_{L^p}),
\la{oneci}
\ee
and summing in $i$ we obtain
\be
\ba
\fr{d}{dt}\sum_{i=1}^m\|q_i\|_{L^p}  \le \C(\|\rho\|_{L^2}^4 +1) \sum_{i=1}^m\|q_i\|_{L^p} + \C m (\|u\|_{L^p} + \|\na\Phi\|_{L^p} +\|\rho\|_{L^p})\\
 \le \C(\|\rho\|_{L^2}^4 +1)\sum_{i=1}^m\|q_i\|_{L^p} + \C m (\|u\|_{L^p} +1)
\ea
\la{l2lp}
\ee
where we used the estimate
\be
\|\na\Phi\|_{L^p}\le \C(\|\rho\|_{L^p} +1) \le \C(\sum_{i=1}^m\|q_i\|_{L^p} +1), 
\la{naphilp}
\ee
from (\ref{normsrhoq}). Thus (\ref{lftylpb}) follows from (\ref{l2lp}), concluding the proof.
\end{proof}
\beg{rem}Because the initial data are in $H^1\subset  L^6$ we have that
(\ref{lftylpb}) holds with $p\in [2,6]$.
\end{rem}
\beg{prop}\la{sufcondR} Let $(c_i, \Phi, u)$ be a strong solution of the system
(\ref{cieq}), (\ref{poi}), (\ref{nse}) (or (\ref{cieq}), (\ref{poi}), (\ref{se})) with boundary conditions (\ref{gammai}), (\ref{phibc}), (\ref{ubc}) on the interval $[0, T]$. Consider the quantity $B(T)$ of (\ref{condB}). Then
\be
\sup_{t\in [0,T]}\sum_{i=1}^m \|c_i(t)\|_{H^1}^2 \le \C\left(\sum_{i=1}^m \|c_i(0)\|_{H^1}^2 +\int_0^T\|u\|_H^2dt+1\right)e^{\C(T+R(T)+U(T))}
\la{h1b}
\ee
holds with
\be
R(T) = \int_0^T\|\rho\|_{L^4}^2dt \le C\left(\sum_{i=1}^m \|c_i(0)\|_{L^4} +\int_0^T\|u\|_{L^4}dt+1\right)^2e^{\C\left(T +B(T)\right)},
\la{RT}
\ee
and
\be
U(T)=\int_0^T\|u\|_{V}^4dt
\la{UT}
\ee
\end{prop}
\begin{proof} In view of (\ref{lftylpb}) with $p=4$ and the fact that $W^{1,p}\subset L^{\infty}$ for $p>3$, we have that
\be
\|\na\Phi\|_{L^{\infty}} \le \C(\|\rho\|_{L^4} +1) \le \C\left(\sum_{i=1}^m \|c_i(0)\|_{L^4}+\int_0^T\|u\|_{L^4}dt +1\right)e^{\C\left(T +B(T)\right)}
\la{naftyphi}
\ee
This is a quantitative bound in terms of the initial data and the constant $B(T)$. Using it together with  
$\|\rho(t)\|_{L^4}$ in the estimate of the evolution of $\|\na q_i\|_{L^2}$ we obtain
\be
\ba
\left|z_iD_i\int_{\Omega}(\na q_i\cdot\na\Phi + q_i\D\Phi)\D q_i dx\right|
\le \C\left (\|\na\Phi\|_{L^{\infty}}\|\na q_i\|_{L^2} + \|q_i\|_{L^4}\|\rho\|_{L^4}\right)\|\D q_i\|_{L^2}\\
\le\C(\|\rho\|_{L^4}+1)\|\na q_i\|_{L^2}\|\D q_i\|_{L^2}
\ea
\la{globnaq}
\ee
instead of (\ref{nonlqh1}), and consequently together with (\ref{conv}), we obtain
\be
\fr{d}{dt}\sum_{i=1}^m\|\na q_i\|^2_{L^2}  + \sum_{i=1}^m \|\D q_i\|_{L^2}^2 \le \C (\|\rho\|_{L^4}^2+\|u\|_V^4+1)\left(\sum_{i=1}^m\|\na q_i\|^2_{L^2}\right)  + \C(\|u\|_H^2+1)
\la{normqnows}
\ee
instead of (\ref{normqs}). Using (\ref{normqnows}) we obtain (\ref{h1b}).
\end{proof}

\beg{thm}\la{bkm} Let $T_1>0$ and let $(c_i, \Phi, u)$ be a strong solution of 
the Nernst-Planck-Stokes system  (\ref{cieq}), (\ref{poi}), (\ref{se}) with boundary conditions (\ref{gammai}), (\ref{phibc}), (\ref{ubc}) on all intervals $[0,T]$ with $T<T_1$. Assume that 
\be
\sup_{T<T_1} \int_0^T\|\rho(t)\|_{L^2}^4dt<\infty.
\la{unifcondB}
\ee
Then there exists $T_2>T_1$ such that $(c_i,\Phi, u)$ can be uniquely continued as a strong solution on $[0,T_2]$.

 Let $T_1>0$ and let $(c_i, \Phi, u)$ be a strong solution of
the Nernst-Planck-Navier-Stokes system  (\ref{cieq}), (\ref{poi}), (\ref{nse}) with boundary conditions (\ref{gammai}), (\ref{phibc}), (\ref{ubc}) on all intervals $[0,T]$ with $T<T_1$. Assume that
\be
\sup_{T<T_1}[\int_0^T\|\rho(t)\|_{L^2}^4dt + \int_0^T\|u\|_V^4dt]<\infty.
\la{unifcondBu}
\ee
Then there exists $T_2>T_1$ such that $(c_i,\Phi, u)$ can be uniquely continued as a strong solution on $[0,T_2]$.
\end{thm}
The converse is obviously also true. Thus (\ref{unifcondB}) (respectively, (\ref{unifcondBu})) is a necessary and sufficient condition for regularity of the Nernst-Planck-Stokes system (respectively, of the Nernst-Planck-Navier-Stokes system).
\beg{proof}
The proof follows directly from Theorem \ref{locex}, Lemma \ref{sufcondB} and Proposition \ref{sufcondR}. We remark that in the case of  (\ref{cieq}), (\ref{poi}), (\ref{se}), $U(T)$ is controlled by $B(T)$.
\end{proof}
\beg{prop}\la{positive} Let $(c_i, \Phi, u)$ be a strong solution of the system
(\ref{cieq}), (\ref{poi}), (\ref{nse}) (or (\ref{cieq}), (\ref{poi}), (\ref{se})) with boundary conditions (\ref{gammai}), (\ref{phibc}), (\ref{ubc}) on the interval $[0, T]$.  If $c_i(x,0)\ge 0$, $i=1, \dots, m$ then $c_i(x,t)\ge 0$ a.e for $t\in [0,T]$.
\end{prop} 
\beg{proof} In order to show this we take a convex function $F:\Rr\to \Rr$ that is nonnegative, twice continuously differentiable, identically zero on the positive semiaxis, and strictly positive on the negative axis. We also assume
\be
F''(y)y^2 \le CF(y)
\la{fass}
\ee
with $C>0$ a fixed constant. Examples of such functions are
\be
F(y) = \left\{
\ba
y^{2m} \quad \quad {\mbox{for}}\quad y<0,\\
0 \quad \quad \quad {\mbox{for}}\quad y\ge 0
\ea
\right.
\la{Fm}
\ee
with $m>1$. (In fact $m=1$ works as well, although we have only $F\in W^{2,\infty}(\Rr)$ in that case.) We multiply the equation (\ref{cieq}) by $F'(c_i)$ and integrate by parts using the fact that $F'(\gamma_i) =0$. We obtain
\be
\fr{d}{dt}\int_{\Omega}F(c_i)dx = -D_i\int_{\Omega} F''(c_i)\left[ |\na c_i|^2  + z_ic_i \na\Phi\cdot\na c_i\right]dx.
\la{intfeq}
\ee
Using a Schwartz inequality and the convexity of $F$, $F''\ge 0$, we have
\be
\fr{d}{dt}\int_{\Omega}F(c_i(x,t))dx \le \fr{CD_i}{2}z_i^2\|\na\Phi\|_{L^{\infty}(\Omega)}^2\int_{\Omega}F(c_i(x,t))dx.
\la{intfineq}
\ee
If $c_i(x,0)\ge 0$ then $F(c_i(x,0))=0$ and (\ref{intfineq}) above shows that 
$F(c_i(x,t))$ has vanishing integral. As $F$ is nonnegative, it follows that $F(c_i(x,t))= 0$ almost everywhere in $x$ and because $F$ does not vanish for negative values it follows that $c_i(x,t)$ is almost everywhere nonnegative.   
\end{proof}
From now on we consider only solutions with $c_i\ge 0$.
\section{Energies}\la{en}
The Navier-Stokes and Stokes energy balance is
\be
\fr{1}{2K}\fr{d}{dt}\int_{\Omega}|u|^2dx  + \fr{\nu}{K}\int_{\Omega}|\na u|^2dx = -\int_{\Omega}\rho(u\cdot\na \Phi)dx.
\la{ens}
\ee
We consider functionals of $(c_i, \Phi)$ which can be used to cancel the right hand side of the Navier-Stokes energy balance.

We denote $(-\D_D)^{-1} $ the inverse of the Laplacian with homogeneous Dirichlet boundary condition. 
We decompose 
\be
\Phi = \Phi_0 + \Phi_W
\la{phis}
\ee
where $\Phi_W$ is harmonic and obeys the inhomogeneous boundary conditions,
\be
\D \Phi_W = 0, \quad {\Phi_W}_{\left | \right. \pa\Omega} = W,
\la{phiw}
\ee
and 
\be 
-\epsilon\D\Phi_0 = \rho, \quad {\Phi_0}_{\left | \right. \pa\Omega} = 0.
\la{phizero}
\ee
so that
\be
\Phi_0 = \fr{1}{\epsilon}(-\D_D)^{-1}\rho.
\la{hpizerorho}
\ee
We introduce
\be
D= \min\{D_1, D_2,\dots D_m\}.
\la{D}
\ee

Let
\be
{\mathcal E}_1 = \int_{\Omega}\left\{\sum_{i=1}^m\Gamma_i\left (\fr{c_i}{\Gamma_i}\log\left(\fr{c_i}{\Gamma_i}\right) - \left(\fr{c_i}{\Gamma_i}\right) +1\right) + \fr{1}{2\epsilon}\rho(-\D_D)^{-1}\rho\right\}dx.
\la{eone}
\ee
\beg{prop}\la{enoneineq}  Let $(c_i, \Phi, u)$ be a strong solution of the system (\ref{cieq}), (\ref{poi}), (\ref{nse}) (or (\ref{cieq}), (\ref{poi}), (\ref{se})) with boundary conditions (\ref{gammai}), (\ref{phibc}), (\ref{ubc}) on the interval $[0, T]$. Then 
\be
\fr{d}{dt}\mathcal E_1 + \mathcal D_1 \le \int_{\Omega}\rho (u\cdot\na \Phi) dx
+\C\left(\sum_{i=1}^m \|c_i-\Gamma_i\|_{L^2} +1\right)(\|u\|_H+1)
\la{dteoneineq}
\ee
holds on $[0,T]$, with
\be
\mathcal D_1 = \fr{D}{2}\int_{\Omega}\left[ \sum_{i=1}^m\left(c_i^{-1}|\na c_i|^2 + z_i^2c_i|\na\Phi|^2\right) + \fr{1}{\epsilon}\rho^2 \right]dx.  
\la{done}
\ee
\end{prop}
The term  $\int_{\Omega}\rho (u\cdot\na\Phi) dx $ in the right hand side of (\ref{dteoneineq}) can be used to cancel the contribution of the electrical forces in the Navier-Stokes energy balance.
\beg{proof} 
We note that
\be
\fr{1}{2\epsilon}\int_{\Omega}\rho(-\D_D)^{-1}\rho dx = \fr{1}{2}\int_{\Omega}\rho \Phi_0 dx. 
\la{poten}
\ee
In order to compute the time evolution of $\mathcal E_1$ we  multiply the equations (\ref{cieq})
by the factors $\log\left(\fr{c_i}{\Gamma_i}\right) + z_i \Phi_0$ and, noting that these factors vanish at the boundary, we integrate by parts:
\be
\ba
\int_{\Omega}((\pa_t + u\cdot\na)c_i)\left(\log\left(\fr{c_i}{\Gamma_i}\right) + z_i \Phi_0\right)dx \\= -D_i\int_{\Omega}c_i\na\left(\log{c_i} + z_i\Phi\right)\cdot \na\left(\log\left(\fr{c_i}{\Gamma_i}\right) + z_i \Phi_0\right)dx \\
= -D_i\int_{\Omega}c_i\left|\na\left(\log{c_i} + z_i\Phi\right)\right|^2dx 
+ D_i\int_{\Omega}c_i\na\left(\log{c_i} + z_i\Phi\right)\cdot\na\left(\log \Gamma_i + z_i\Phi_W\right) dx
\ea
\la{eoneone}
\ee
We have thus
\be
\ba
\int_{\Omega}((\pa_t + u\cdot\na)c_i)\left(\log\left(\fr{c_i}{\Gamma_i}\right) + z_i \Phi_0\right)dx \\\le -\fr{1}{2}D_i\int_{\Omega}c_i\left|\na\left(\log{c_i} + z_i\Phi\right)\right|^2dx + 
\fr{1}{2} D_i\int_{\Omega}c_i|\na(\log \Gamma_i + z_i \Phi_W)|^2 dx
\ea
\la{eoneoneineone}
\ee
In view of the fact that
\be
((\pa_t + u\cdot\na)c_i) \log\left(\fr{c_i}{\Gamma_i}\right) = 
(\pa_t + u\cdot\na)\left(c_i \log\left(\fr{c_i}{\Gamma_i}\right) -c_i\right) +
c_i u\cdot\na\log\Gamma_i,
\la{fact}
\ee
summing in $i$, on  the left hand side we have
\be
\ba
\sum_{i=1}^m\int_{\Omega}((\pa_t + u\cdot\na)c_i)\left(\log\left(\fr{c_i}{\Gamma_i}\right) + z_i \Phi_0\right)dx \\ 
=\fr{d}{dt}\int_{\Omega}\sum_{i=1}^mc_i(\log\left(\fr{c_i}{\Gamma_i}\right) -1) dx + \int_{\Omega}(\pa_t\sum_{i=1}^m(z_ic_i))\Phi_0 dx\\ +\int_{\Omega}(u\cdot\na (\sum_{i=1}^m z_ic_i))\Phi_0 dx  + \int_{\Omega}\sum_{i=1}^2 c_i u\cdot\na \log\Gamma_i dx\\
=\fr{d}{dt}\mathcal E_1 + \int_{\Omega}(u\cdot\na\rho) \Phi_0dx + \sum_{i=1}^m
\int_{\Omega}c_i u\cdot\na\log\Gamma_i dx.
\ea
\la{leftone}
\ee
In the last equality we used
\be
\fr{d}{dt}\fr{1}{2\epsilon}\int_{\Omega}\rho(-\D_D)^{-1}\rho dx = \int_{\Omega}(\pa_t\rho)\Phi_0 dx
\la{dtpoten}
\ee
because $(-\D_D)^{-1}$ is selfadjoint. Combining (\ref{eoneoneineone}) and (\ref{leftone}) we obtain
\be
\ba
\fr{d}{dt}\mathcal E_1  \le  -\fr{1}{2}\sum_{i=1}^mD_i\int_{\Omega}c_i\left|\na\left(\log{c_i} + z_i\Phi\right)\right|^2dx + \fr{1}{2}\sum_{i=1}^m D_i\int_{\Omega}c_i|\na(\log \Gamma_i + z_i\Phi_W)|^2 dx \\ - \sum_{i=1}^m\int_{\Omega}c_i u\cdot\na\log\Gamma_i dx
- \int_{\Omega}(u\cdot\na\rho)\Phi_0dx\\
= -\fr{1}{2}\sum_{i=1}^mD_i\int_{\Omega}c_i\left|\na\left(\log{c_i} + z_i\Phi\right)\right|^2dx + \fr{1}{2}\sum_{i=1}^m D_i\int_{\Omega}c_i|\na(\log \Gamma_i + z_i \Phi_W)|^2 dx \\ - \sum_{i=1}^m\int_{\Omega}c_i u\cdot\na\log\Gamma_i dx + \int_{\Omega}\rho u\cdot\na (\Phi -\Phi_W)dx 
\ea
\la{enoneone}
\ee
We note that
\be
\ba
\fr{1}{2}\sum_{i=1}^mD_i\int_{\Omega}c_i\left|\na\left(\log{c_i} + z_i\Phi\right)\right|^2dx\ge \fr{D}{2}\sum_{i=1}^m\int_{\Omega}c_i\left|\na\left(\log{c_i} + z_i\Phi\right)\right|^2dx\\
=\fr{D}{2}\sum_{i=1}^m\int_{\Omega}(c_i^{-1}|\na c_i|^2 + z_i^2c_i|\na\Phi|^2)dx + D\int_{\Omega}\sum_{i=1}^m z_i\na c_i\cdot \na \Phi dx \\
=\fr{D}{2}\sum_{i=1}^m\int_{\Omega}(c_i^{-1}|\na c_i|^2 + z_i^2c_i|\na\Phi|^2)dx + D\int_{\Omega}\sum_{i=1}^m z_i\na (c_i-\Gamma_i)\cdot \na \Phi dx   \\
+ D\int_{\Omega}\sum_{i=1}^m z_i\Gamma_i\cdot \na \Phi dx \\
=\fr{D}{2}\sum_{i=1}^m\int_{\Omega}(c_i^{-1}|\na c_i|^2 + z_i^2c_i|\na\Phi|^2)dx + D\int_{\Omega}\sum_{i=1}^m z_i(c_i-\Gamma_i)( -\D\Phi) dx   \\
+ D\int_{\Omega}\sum_{i=1}^m z_i\Gamma_i\cdot \na \Phi dx \\
=\fr{D}{2}\sum_{i=1}^m\int_{\Omega}(c_i^{-1}|\na c_i|^2 + z_i^2c_i|\na\Phi|^2)dx + D\fr{1}{\epsilon}\int_{\Omega}\sum_{i=1}^m z_i(c_i-\Gamma_i)\rho dx   \\
+ D\int_{\Omega}\sum_{i=1}^m z_i\Gamma_i\cdot \na \Phi dx \\
=\fr{D}{2}\sum_{i=1}^m\int_{\Omega}(c_i^{-1}|\na c_i|^2 + z_i^2c_i|\na\Phi|^2)dx + D\fr{1}{\epsilon}\int_{\Omega}(\rho -\sum_{i=1}^m z_i\Gamma_i)\rho dx   \\
+ D\int_{\Omega}\sum_{i=1}^m z_i\Gamma_i\cdot \na \Phi dx \\
\ge \fr{D}{2}\sum_{i=1}^m\int_{\Omega}(c_i^{-1}|\na c_i|^2 + z_i^2c_i|\na\Phi|^2)dx + \fr{D}{2\epsilon}\int_{\Omega}\rho^2 dx\\
-\fr{D}{2\epsilon}\int_{\Omega}\left |\sum_{i=1}^mz_i\Gamma_i\right|^2dx
 + D\int_{\Omega}\sum_{i=1}^m z_i\Gamma_i\cdot \na \Phi dx.
\ea
\la{dissiplow}
\ee
From (\ref{enoneone}) and (\ref{dissiplow}) we obtain
\be
\fr{d}{dt}\mathcal E_1 + \mathcal D_1 \le \mathcal Q_1 + \int_{\Omega}\rho u\cdot\na\Phi dx 
\la{eonebalance}
\ee
with $\mathcal D_1$ given in (\ref{done})
and 
\be
\ba
\mathcal Q_1 =  \fr{1}{2}\sum_{i=1}^m D_i\int_{\Omega}c_i|\na(\log \Gamma_i + z_i\Phi_W)|^2 dx\\ - \sum_{i=1}^m\int_{\Omega}c_i u\cdot\na\log\Gamma_i dx - \int_{\Omega}\rho u\cdot\na \Phi_W dx\\ 
+\fr{D}{2\epsilon}\int_{\Omega}\left | \sum_{i=1}^mz_i\Gamma_i\right|^2dx
 - D\int_{\Omega}\sum_{i=1}^m z_i\Gamma_i\cdot \na \Phi dx.
\ea
\la{qone}
\ee
Note that $\mathcal Q_1$ is at most quadratic in terms of the unknowns $u, c_i$, in view of the fact that both $\rho$ and $\Phi$ are affine in $c_i$. The inequality (\ref{dteoneineq}) follows by bounding $Q_1$. \end{proof}

A useful energy is the potential part in $\mathcal E_1$:
\be
\mathcal P = \fr{1}{2\epsilon} \int_{\Omega}\rho (-\D_D)^{-1}\rho dx.
\la{P}
\ee

\beg{prop}\la{potentialene}  Let $(c_i, \Phi, u)$ be a strong solution of the system (\ref{cieq}), (\ref{poi}), (\ref{nse}) (or (\ref{cieq}), (\ref{poi}), (\ref{se})) with boundary conditions (\ref{gammai}), (\ref{phibc}), (\ref{ubc}) on the interval $[0, T]$. Then
\be
\fr{d}{dt}\mathcal P + \mathcal D_2 \le  \int_{\Omega}\rho (u\cdot\na\Phi) dx
+ \C\left(\sum_{i=1}^m\|c_i-\Gamma_i\|_{L^2}+1\right)(\|\rho\|_{L^2}+1) +
\C\|\rho\|_{L^2}\|u\|_H
\la{pineq}
\ee
holds on $[0,T]$, with
\be
\mathcal D_2 =  \fr{1}{2}\sum_{i=1}^m  z_i^2D_i\int_{\Omega}c_i \left| \na\Phi\right|^2dx.
\la{d2}
\ee
\end{prop}
We note that the term $\int_{\Omega}\rho(u\cdot\na \Phi)dx$ in the right hand side of (\ref{pineq}) can  be used to cancel the contribution of electrical forces in the Navier-Stokes energy balance. 

\beg{proof}
In order to compute the time evolution of (\ref{P}) we take the equations (\ref{cieq}), multiply by the factors $z_i\Phi_0$ and  integrate by parts in view of the fact that $\Phi_0$ vanishes on the boundary. We obtain
\be
\ba
\int_{\Omega}((\pa_t + u\cdot\na)c_i)z_i\Phi_0dx = -D_iz_i\int_{\Omega}(\na c_i+ z_i
 c_i\na\Phi)\cdot\na\Phi_0dx\\ = -D_iz_i\int_{\Omega}\na c_i\cdot\na\Phi_0 -z_i^2D_i
\int_{\Omega}c_i\na\Phi\cdot\na\Phi_0 dx\\
= -z_iD_i\int_{\Omega}\na(c_i-\Gamma_i)\cdot\na \Phi_0 dx -D_iz_i\int_{\Omega}\na\Gamma_i\cdot\na\Phi_0\\ - z_i^2D_i\int_{\Omega}c_i \left| \na\Phi\right|^2dx + z_i^2D_
i\int_{\Omega}c_i  \na\Phi \cdot\na\Phi_Wdx\\
=-D_iz_i\epsilon^{-1}\int_{\Omega}(c_i-\Gamma_i)\rho dx   - z_i^2D_i\int_{\Omega}c_i
 \left| \na\Phi\right|^2dx\\ + z_i^2D_i\int_{\Omega}c_i  \na\Phi \cdot\na\Phi_Wdx -D
_iz_i\int_{\Omega}\na\Gamma_i\cdot\na\Phi_0dx
\ea
\la{potrightone}
\ee
In the last equality we used the fact that $c_i-\Gamma_i$ vanishes on the boundary and the fact that $-\epsilon\D \Phi_0 = \rho$. Summing in $i$, on the left  hand side
 we have
\be
\ba
\sum_{i=1}^m\int_{\Omega}((\pa_t + u\cdot\na)c_i)z_i\Phi_0dx =
\int_{\Omega}(\pa_t\rho)\Phi_0dx  + \int_{\Omega}(u\cdot\na\rho)\Phi_0dx\\
=\fr{d}{dt}\mathcal P +\int_{\Omega}(u\cdot\na\rho)\Phi_0 dx
\ea
\la{leftpone}
\ee
Putting together (\ref{potrightone}) and (\ref{leftpone}) 
\be
\ba
\fr{d}{dt}\mathcal P + \sum_{i=1}^m  z_i^2D_i\int_{\Omega}c_i \left| \na\Phi\right|^2dx\\
= -\sum_{i=1}^mD_iz_i\epsilon^{-1}\int_{\Omega}(c_i-\Gamma_i)\rho dx\\
+\sum_{i=1}^m z_i^2D_i\int_{\Omega}c_i  \na\Phi \cdot\na\Phi_Wdx -
\sum_{i=1}^mD_iz_i\int_{\Omega}\na\Gamma_i\cdot\na\Phi_0dx\\
-\int_{\Omega}\rho(u\cdot\na \Phi_W) dx  
+ \int_{\Omega}\rho u\cdot\na\Phi dx.
\ea
\la{pbalanceone}
\ee
After a Schwartz inequality we obtain
\be
\fr{d}{dt}\mathcal P + \mathcal D_2 \le  \mathcal Q_{2} + \int_{\Omega}\rho u\cdot\na\Phi dx
\la{pbalance}
\ee
where $\mathcal D_2$ is given in (\ref{d2}) and
\be
\ba
\mathcal Q_{2} = -\sum_{i=1}^mD_iz_i\epsilon^{-1}\int_{\Omega}(c_i-\Gamma_i)\rho dx\\
+\fr{1}{2}\sum_{i=1}^m z_i^2D_i\int_{\Omega}c_i \left| \na\Phi_W\right |^2dx -
\sum_{i=1}^mD_iz_i\int_{\Omega}\na\Gamma_i\cdot\na\Phi_0dx -\int_{\Omega}\rho(u\cdot\na \Phi_W) dx  
\ea
\la{q2}
\ee
Unlike the term $\mathcal Q_1$ of (\ref{qone}), $\mathcal Q_2$ has no $(u,c)$ quadratic terms, the only quadratic terms are of the type $(c,\rho)$ or $(u,\rho)$ (the  $(u,\Phi_0)$ term is is of $(u,\rho)$ type in this accounting).
Estimating $Q_2$ we obtain (\ref{pineq}).

\end{proof}

\section{Quadratic bounds}\la{glob}
We estimate the sum of $L^2$ norms of $c_i.$ We take the scalar product of the equations (\ref{Fiqieq}) with $\fr{1}{D_i} q_i$ and add. We obtain first
\be
\fr{d}{dt}\sum_{i=1}^m  \fr{1}{2D_i}\int_{\Omega} q_i^2 dx + \sum_{i=1}^m \int_{\Omega}|\na q_i|^2dx = -\fr{1}{2}\sum_{i=1}^m  z_i\int_{\Omega}\na\Phi\cdot\na(q_i^2)dx + \sum_{i=1}^m  \fr{1}{D_i}\int_{\Omega}F_iq_idx.  
\la{enstep1}
\ee 
The integartion by parts is justified because of (\ref{qibc}). We integrate by parts one more time using the same boundary conditions and (\ref{poi})
\be
\fr{d}{dt}\sum_{i=1}^m  \fr{1}{2D_i}\int_{\Omega} q_i^2 dx + \sum_{i=1}^m \int_{\Omega}|\na q_i|^2dx = -\fr{1}{2\epsilon}\int_{\Omega}\rho\sum_{i=1}^m  z_i q_i^2 dx + \sum_{i=1}^m  \fr{1}{D_i}\int_{\Omega}F_iq_idx.
\la{enstep2}
\ee

\beg{thm}\la{twospecies} Consider  $m=2$, $z_1=1$, $z_2=-1$. Let $T>0$ be arbitrary. Let $c_i(\cdot, 0)>0$, $c_i(\cdot, 0)\in H^1$, ${c_i}_{|\pa\Omega} = \gamma_i$ and $u_0\in V$ be given. Then the system  (\ref{cieq}), (\ref{poi}), (\ref{se})) with boundary conditions (\ref{gammai}), (\ref{phibc}), (\ref{ubc}) has global strong solutions on $[0,T]$.
The system (\ref{cieq}), (\ref{poi}), (\ref{nse}) has global strong solutions if
\be
\int_0^T\|u\|_V^4dt <\infty.
\la{nsesuf} 
\ee
Moreover 
\be
\sup_{0\le t \le T}\sum_{i=1}^2 \|c_i(t)\|_{H^1}^2 + \int_0^T\sum_{i=1}^2\|c_i(t)\|^2_{H^2}dt \le \C\left [\sum_{i=1}^2 \|c_i(0)\|_{H^1}^2+\int_0^T\|u\|_H^2dt +1\right]e^{\C (T+R(T)+U(T))}
\la{strongbm2}
\ee
holds for all $T$, where $R(T)$ is given by (\ref{RT}), $U(T)$ is given by (\ref{UT}), and  with $\C$ depending only on the boundary conditions $\gamma_i$ and $W$, domain $\Omega$, and parameters $\nu, D_i, \epsilon, K$.
\end{thm}
\beg{proof}
We note that, when $m=2$ and $z_1 =1$, $z_2=-1$, then
\be
\sum_{i=1}^2  z_i q_i^2 = (\rho  -\Gamma_1+\Gamma_2 )(c_1+ c_2 -\Gamma_1 -\Gamma_2). 
\la{rhosum}
\ee
Thus 
\be
\ba
\rho\sum_{i=1}^2  z_i q_i^2 = \rho^2(c_1+c_2)  - (\Gamma_1-\Gamma_2)\rho(q_1+q_2) - \rho^2 (\Gamma_1 + \Gamma_2)\\
\ge |\rho|^3 -(\Gamma_1 - \Gamma_2)\rho(q_1+q_2) -\rho^2(\Gamma_1 +\Gamma_2) 
\ea
\la{lowcube}
\ee
because 
\be
c_1+ c_2 \ge |\rho|.
\la{sumdiff}
\ee
Now we use H\"{o}lder and Young inequalities to bound in (\ref{enstep2})
\be
\ba
\fr{1}{2\epsilon}\int_{\Omega}\rho\sum_{i=1}^2  z_i q_i^2 dx \ge
\fr{1}{2\epsilon}\int_{\Omega}|\rho|^3dx\\ - \fr{1}{2\epsilon}(\|q_1\|_{L^2(\Omega)}+\|q_2\|_{L^2(\Omega)})\|\Gamma_2-\Gamma_1\|_{L^6(\Omega)}\|\rho\|_{L^3(\Omega)} \\
-\fr{1}{2\epsilon}[\|\Gamma_1\|_{L^3(\Omega)} + \|\Gamma_2\|_{L^3(\Omega)}]\|\rho\|_{L^3(\Omega)}^2\\
\ge \fr{1}{4\epsilon}\|\rho\|_{L^3}^3 - \fr{1}{4 L^2}\|q_1\|_{L^2}^2 -\fr{1}{4 L^2}\|q_2\|_{L^2}^2 -\C  \\
\ea
\la{inter}
\ee
with $L$ the constant in the Poincar\'{e} inequality
\be
\|\na q\|_{L^2(\Omega)}^2 \ge L^{-2}\|q\|^2_{L^2(\Omega)}.
\la{poin}
\ee
From (\ref{enstep2}), (\ref{inter}) and (\ref{poin}) we obtain
\be
\ba
\fr{d}{dt}\sum_{i=1}^2  \fr{1}{2D_i}\int_{\Omega} q_i^2 dx + \fr{3}{4}\sum_{i=1}^2 \int_{\Omega}|\na q_i|^2dx +  \fr{1}{4\epsilon}\|\rho\|_{L^3}^3\\
\le \C + \sum_{i=1}^2\fr{1}{D_i}\int_{\Omega}F_iq_idx
\ea
\la{intern}
\ee
We have
\be
\left|\int_{\Omega}F_iq_i dx\right| \le \C(\|\rho\|_{L^2} +\|u\|_H +1)\|q_i\|_{L^2}
\la{fqbound}
\ee
and therefore we have that
\be
\mathcal E_3 = \sum_{i=1}^2  \fr{1}{D_i}\int_{\Omega} q_i^2 dx
\la{e3}
\ee
obeys
\be
\fr{d}{dt}\mathcal E_3 + \mathcal D_3 \le  \C + \widetilde{\C}\|u\|_H^2
\la{e3ineq}
\ee
with
\be
\mathcal D_3 =  \sum_{i=1}^2\fr{1}{2} \int_{\Omega}|\na q_i|^2dx +  \fr{1}{4\epsilon}\|\rho\|_{L^3}^3.
\la{d3}
\ee
We singled out the coefficient $\widetilde{\C}$ of $\|u\|_H^2$ because we use it next. We take a constant
\be
\delta = \fr{\nu}{2KL^2\widetilde{\C}}
\la{deltaLC}
\ee 
such that the dissipation in the Navier-Stokes energy balance exceeds twice the contribution from $\|u\|_H^2$ in the right hand side of (\ref{e3ineq}) when the latter is multiplied by $\delta$, 
\be
\fr{\nu}{K}\int_{\Omega}|\na u|^2dx \ge 2\delta\widetilde{\C}\|u\|_H^2.
\la{dissbeats}
\ee
We consider
\be
\mathcal F = \fr{1}{2K}\|u\|_H^2 + \mathcal P + \delta \mathcal E_3
\la{mathcalf}
\ee
and, using (\ref{ens}), (\ref{pineq}) and (\ref{e3ineq}) multiplied by $\delta$ we obtain
\be
\ba
\fr{d}{dt}\mathcal F + \fr{\nu}{2K}\int_{\Omega}|\na u|^2dx + 
\fr{\delta}{2}\sum_{i=1}^2\|\na q_i\|^2_{L^2} + \fr{\delta}{4\epsilon}\|\rho\|_{L^3}^3\\
 \le \C[\|\rho\|_{L^2}(\sum_{i=1}^2\|q_i\|_{L^2} + 1) + \|\rho\|_{L^2}\|u\|_H +  \sum_{i=1}^2\|q_i\|_{L^2} + 1].
\ea
\la{mathcalfineq}
\ee
The positive cubic term in $\rho$ on the left hand side together with the rest of positive quadratic dissipative terms on the left hand side can be used to absorb all the  quadratic terms on the right hand side, because they all involve at least one $\rho$, and the linear terms are also absorbed using  Poincar\'{e} inequalities for both $q_i$ and for $u$. This results in
\be
\fr{d}{dt}\mathcal F + c_{\Gamma} \mathcal F \le \C
\la{mathcalfinal}
\ee
with $c_{\Gamma}>0$. It follows that
\be
\mathcal F(t) \le \mathcal F(0) e^{-c_{\Gamma}t} + \C
\la{mathcalFbound}
\ee
This implies in particular that
\be
\|\rho(t)\|_{L^2}^2\le \C\mathcal F(0)e^{-c_{\Gamma}t} + \C
\la{rhol2b}
\ee
and, integrating in time, (\ref{condB}) holds
\be
B(T)\le \C(\mathcal F(0)+ T).
\la{BTB}
\ee
Moreover, the dissipation is time integrable, 
\be
\int_0^T\left\{ \fr{\nu}{2K}\int_{\Omega}|\na u|^2dx +
\fr{\delta}{2}\sum_{i=1}^2\|\na q_i\|^2_{L^2} + \fr{\delta}{4\epsilon}\|\rho\|_{L^3}^3\right\}dt\le \C(\mathcal F(0) +T). 
\la{dissip}
\ee
It follows from (\ref{h1b}) that (\ref{strongbm2}) holds.
\end{proof}

\beg{thm}\la{equaldiffs} Consider  $z_i=\pm 1$, $i=1, \dots,m$, and assume $D_1 = D_2 =\dots = D_m =D>0$. Let $T>0$ be arbitrary. Let $c_i(\cdot, 0)>0$, $c_i(\cdot, 0)\in H^1$, ${c_i}_{|\pa\Omega} = \gamma_i$ and $u_0\in V$ be given. Then the system  (\ref{cieq}), (\ref{poi}), (\ref{se})) with boundary conditions (\ref{gammai}), (\ref{phibc}), (\ref{ubc}) has global strong solutions on $[0,T]$. The system (\ref{cieq}), (\ref{poi}), (\ref{nse}) has global stromg solutions if
\be
\int_0^T\|u\|_V^4dt <\infty.
\la{nsesuff}
\ee
Moreover
\be
\sup_{0\le t \le T}\sum_{i=1}^2 \|c_i(t)\|_{H^1}^2 + \int_0^T\sum_{i=1}^2\|c_i(t)\|^2_{H^2}dt \le \C\left [\sum_{i=1}^2 \|c_i(0)\|_{H^1}^2 +\int_0^T\|u\|_H^2dt+1\right]e^{\C (T+R(T)+U(T))}
\la{strongbeqdiff}
\ee
holds for all $T$, where $R(T)$ is given by (\ref{RT}), $U(T)$ is given by (\ref{UT}), and with $\C$ depending only on the boundary conditions $\gamma_i$ and $W$, domain $\Omega$, and parameters $\nu, D_i, \epsilon, K$.
\end{thm}
\beg{proof}
We consider the auxiliary variables 
\be
S = \sum_{i=1}^M q_i
\la{S}
\ee
and 
\be
Z = \sum_{i=1}^mz_iq_i.
\la{Z}
\ee
Summing in (\ref{Fiqieq}) we have 
\be
\left\{
\ba
\left(\pa_t +u\cdot\na\right) S = D\left(\D S +\div(Z\na\Phi)\right) + F_S\\
\left(\pa_t + u\cdot\na\right)Z = D(\left(\D Z + \div(S\na\Phi)\right) + F_Z,
\ea
\right.
\la{SZeq}
\ee
with 
\be
F_S =  \sum_{i=1}^m F_i,
\la{FS}
\ee
\be
F_Z = \sum_{i=1}^mz_iF_i,
\la{FZ}
\ee
and $F_i$ given in (\ref{Fi}) and with $D$ in (\ref{D}). Multiplying by $S$ and $Z$
and integrating by parts (twice in the nonlinear term, once in linear terms) we obtain
\be
\fr{1}{2}\fr{d}{dt}\int_{\Omega}(S^2 + Z^2)dx + D\int_{\Omega}(|\na S|^2 + |\na Z|^2)dx
 = -\fr{D}{\epsilon}\int_{\Omega} SZ\rho dx + \int_{\Omega}\left(SF_S + ZF_Z\right)dx
\la{enSZ}
\ee
Now we use
\be
Z = \rho  - \Gamma_Z
\la{Zrho}
\ee
and 
\be
S = \sum_{i=1}^M c_i - \Gamma_S
\la{Ssum}
\ee
with
\be
\Gamma_S = \sum_{i=1}^M \Gamma_i,
\la{gammas}
\ee
and 
\be
\Gamma_Z = \sum_{i=1}^m z_i\Gamma_i,
\la{gammaz}
\ee
together with
\be
\sum_{i=1}^m c_i \ge |\rho|,
\la{Srhon}
\ee
to deduce
\be
SZ\rho \ge |\rho|^3 - \Gamma_S \rho^2 - S\rho \Gamma_Z.
\la{szrho}
\ee
Let us note the relationships
\be
\left\{
\ba
F_S = - u\cdot\na \Gamma_S  + D(\D \Gamma_S + \div(\Gamma_Z\na\Phi)),\\
F_Z = - u\cdot\na \Gamma_Z   + D(\D \Gamma_Z + \div(\Gamma_S\na\Phi)).\\
\ea
\right.
\la{FGamma}
\ee
We deduce that
\be
\ba
\left| \int_{\Omega}(SF_S +ZF_Z)dx\right | \le\fr{1}{2}\int_{\Omega}(|\na S|^2 + |\na Z|^2)dx \\ + \int_{\Omega}\left[|u|^2\left(\Gamma_S^2 + \Gamma_Z^2\right) + D^2|\na \Gamma_S|^2 + D^2|\na\Gamma_Z|^2\right ]dx + D^2\int_{\Omega}\left(|\Gamma_S|^2 + |\Gamma_Z|^2\right)|\na\Phi|^2dx
\ea
\la{sfzf}
\ee
Using(\ref{enSZ}), (\ref{szrho}), (\ref{sfzf}) and  (\ref{naphilp}) we obtain
\be
\ba
\fr{d}{dt}\int_{\Omega}(S^2 + Z^2)dx + \fr{D}{2}\int_{\Omega}(|\na S|^2 + |\na Z|^2)dx
 +\fr{D}{4\epsilon}\int_{\Omega} |\rho|^3 dx \\
\le \C + \widetilde{\C}\|u\|_H^2
\ea
\la{ensszb}
\ee
We take $\delta$ defined in (\ref{deltaLC}) with the current $\widetilde{\C}$ and consider the functional 
\be
\mathcal G = \fr{1}{2K}\|u\|_H^2 + \mathcal P + \delta\int_{\Omega}(S^2 + Z^2)dx
\la{mathcalG}
\ee
and obtain from  (\ref{ens}), (\ref{pineq}) and (\ref{ensszb})
\be
\ba
\fr{d}{dt}\mathcal G + \fr{\nu}{2K}\int_{\Omega}|\na u|^2dx +
\fr{\delta}{2}(\|\na S\|^2_{L^2} + \|\na Z\|_{L^2}^2)+  \fr{\delta}{4\epsilon}\|\rho\|_{L^3}^3\\
 \le \C[\|\rho\|_{L^2}(\sum_{i=1}^2\|c_i-\Gamma_i\|_{L^2} + 1) + \|\rho\|_{L^2}\|u\|_H +  \sum_{i=1}^2\|c_i-\Gamma_i\|_{L^2} + 1].
\ea
\la{gineq}
\ee
Now we note that
\be
0\le c_i \le \sum_{i=1}^m c_i = S + \Gamma_S
\la{cs}
\ee
implies that
\be
\|c_i-\Gamma_i\|_{L^2} \le \|S\|_{L^2} + \C
\la{qiS}
\ee
and the Poincare inequality for $S$ implies
\be
\|\na S\|_{L^2}^2 \ge \fr{1}{L^2}\|c_i-\Gamma_i\|^2_{L^2} -\C.
\la{pois}
\ee
Therefore we obtain using Poincar\'{e}, H\"{o}lder and Young inequalities,
\be
\fr{d}{dt}\mathcal G + c_{\Gamma}\mathcal G \le \C
\la{ginedecay}
\ee
with $c_{\Gamma}>0$. The rest of the proof follows as in the proof of Theorem \ref{twospecies}.
\end{proof}

\vspace{.5cm}

{\bf{Acknowledgment.}} The work of PC was partially supported by NSF grant DMS-
171398.

\end{document}